\documentclass[11pt]{article}
\usepackage{amssymb,amsmath,latexsym}
\usepackage{pdfsync,color}

\allowdisplaybreaks

\begin{document}

\begin{doublespace}

\def\1{{\bf 1}}
\def\ind{{\bf 1}}
\def\nn{\nonumber}

\def\sA {{\cal A}} \def\sB {{\cal B}} \def\sC {{\cal C}}
\def\sD {{\cal D}} \def\sE {{\cal E}} \def\sF {{\cal F}}
\def\sG {{\cal G}} \def\sH {{\cal H}} \def\sI {{\cal I}}
\def\sJ {{\cal J}} \def\sK {{\cal K}} \def\sL {{\cal L}}
\def\sM {{\cal M}} \def\sN {{\cal N}} \def\sO {{\cal O}}
\def\sP {{\cal P}} \def\sQ {{\cal Q}} \def\sR {{\cal R}}
\def\sS {{\cal S}} \def\sT {{\cal T}} \def\sU {{\cal U}}
\def\sV {{\cal V}} \def\sW {{\cal W}} \def\sX {{\cal X}}
\def\sY {{\cal Y}} \def\sZ {{\cal Z}}

\def\bA {{\mathbb A}} \def\bB {{\mathbb B}} \def\bC {{\mathbb C}}
\def\bD {{\mathbb D}} \def\bE {{\mathbb E}} \def\bF {{\mathbb F}}
\def\bG {{\mathbb G}} \def\bH {{\mathbb H}} \def\bI {{\mathbb I}}
\def\bJ {{\mathbb J}} \def\bK {{\mathbb K}} \def\bL {{\mathbb L}}
\def\bM {{\mathbb M}} \def\bN {{\mathbb N}} \def\bO {{\mathbb O}}
\def\bP {{\mathbb P}} \def\bQ {{\mathbb Q}} \def\bR {{\mathbb R}}
\def\bS {{\mathbb S}} \def\bT {{\mathbb T}} \def\bU {{\mathbb U}}
\def\bV {{\mathbb V}} \def\bW {{\mathbb W}} \def\bX {{\mathbb X}}
\def\bY {{\mathbb Y}} \def\bZ {{\mathbb Z}}
\def\R {{\mathbb R}} \def\RR {{\mathbb R}}
\def\n{{\bf n}}

\newcommand{\expr}[1]{\left( #1 \right)}
\newcommand{\cl}[1]{\overline{#1}}
\newtheorem{thm}{Theorem}[section]
\newtheorem{lemma}[thm]{Lemma}
\newtheorem{defn}[thm]{Definition}
\newtheorem{prop}[thm]{Proposition}
\newtheorem{corollary}[thm]{Corollary}
\newtheorem{remark}[thm]{Remark}
\newtheorem{example}[thm]{Example}
\newtheorem{assumption}[thm]{Assumption}
\numberwithin{equation}{section}
\def\ee{\varepsilon}
\def\qed{{\hfill $\Box$ \bigskip}}
\def\NN{{\cal N}}
\def\AA{{\cal A}}
\def\MM{{\cal M}}
\def\BB{{\cal B}}
\def\CC{{\cal C}}
\def\LL{{\cal L}}
\def\DD{{\cal D}}
\def\FF{{\cal F}}
\def\EE{{\cal E}}
\def\QQ{{\cal Q}}
\def\RR{{\mathbb R}}
\def\R{{\mathbb R}}
\def\L{{\bf L}}
\def\K{{\bf K}}
\def\S{{\bf S}}
\def\A{{\bf A}}
\def\E{{\mathbb E}}
\def\F{{\bf F}}
\def\P{{\mathbb P}}
\def\N{{\mathbb N}}
\def\eps{\varepsilon}
\def\wh{\widehat}
\def\wt{\widetilde}
\def\pf{\noindent{\bf Proof.} }
\def\pff{\noindent{\bf Proof} }
\def\beq{\begin{equation}}
\def\eeq{\end{equation}}
\def\bee{\begin{equation}}
\def\eee{\end{equation}}
\def\osc{\mathrm{Osc}}

\title{\Large \bf Regularity for fully nonlinear equations driven by spatial-inhomogeneous nonlocal operators}

\author{
{\bf Jongchun Bae}\thanks{This work was supported by the National Research Foundation of Korea(NRF) grant funded by the Korea government(MEST) (2013004822)} }

\date{}
\maketitle

\begin{abstract}
In this paper we consider a large class of fully nonlinear integro-differential equations. The class of our nonlocal operators we consider is not spatial homogeneous and we put   mild assumptions on its kernel near zero.  We prove the H\"older regularity for such equation. In particular, our result covers the case that the kernel $K(x,y)$ is comparable to $|x-y|^{-d-\alpha} \ln (|x-y|^{-1})$ for $|x-y|<c$ where $0<\alpha<2$.
\end{abstract}

\noindent {\bf AMS 2010 Mathematics Subject Classification}: Primary  60J75, 47G20, 
Secondary  35B65.

\noindent {\bf Keywords and phrases:}
L\'evy processes, nonlocal operators, fully nonlinear equation, H\"older continuity, integro-differential equations

\section{Introduction}\label{S:intro}
For many models describing scientific phenomena it is important to find  whether a solution of model equation is continuous or not. The regularity theory for second-order elliptic partial differential operators in divergence form was developed in mid-20th century by many mathematicians: Morrey, De Giorgi, Nash, Moser and so on. For non-divergence type operators, Krylov and Safonov \cite{KrSa1, KrSa2}
proved the corresponding results
 using probabilistic methods. A non-divergence type differential  operator may be represented by infinitesimal generator of a diffusion process which is a  strong Markov process with continuous sample paths. The simplest example is the Laplace operator which is the infinitesimal generator of Brownian motion.

Since researchers in the areas of applied mathematics  found that real world phenomena are frequently described well if one considers jump processes,
many scientists have been interested in integro-differential operators which can be regarded as the infinitesimal generator of  processes with discontinuous sample paths. Analytically, in \cite{S06} Silvestre proved the H\"older estimate for solutions in proper sense of integro-differential equations with a kernel comparable to that of the fractional Laplacian. Kassmann obtained the same result for weak solution by developing a nonlocal version of De Giorgi-Nash-Moser theorem in \cite{K09}. Caffarelli and Silvestre established regularity theory for fully nonlinear integro-differential equations extending the results for elliptic partial differential equations in \cite{CS09}. Overcoming the difficulty of non-symmetry of the kernels, the authors in \cite{KL12, KL13, CD12} obtained  
H\"older and $C^{1,\alpha}$ estimates for the case of non-symmetric kernels.  All of them basically considered families of kernels comparable to that of the fractional Laplacian. 
In this paper we will consider more general kernels than the kernel of the fractional Laplacian.

In probabilistic point of view it is important to estimate transition probabilities and jump measures of Markov 
processes with discontinuous sample paths. 
Such estimates are closely related to Harnack inequality and H\"older estimate for harmonic functions with respect to the  processes, which have been active research areas over last decade and more. In \cite{BK05a, BK05b, BL02, SV04} the Harnack inequality and the H\"older estimate of nonnegative bounded harmonic functions were obtained by probabilistic methods. 
In \cite{CKK2, CKK3, CK} H\"older estimates and the two-sided estimates of transition densities  of  general Markov processes
were obtained.

Let us now  set up some notations.
We use ``$:=$" to denote a definition, which
is  read as ``is defined to be"; we denote $a \wedge b := \min \{ a, b\}$;
we  use $dx$ to denote the Lebesgue measure in $\bR^d$; for a Borel set $A\subset \bR^d$, we use $|A|$ to denote its Lebesgue measure and $\1_{A}$ to denote its indicator function;
 we 
denote by $B(x, r)$ the open
ball centered at $x\in \bR^d$ with radius $r>0$; for every function $f$, let
$f_+:=f \vee 0$.

We consider a measurable function $J : \bR^d \times \bR^d\setminus \{(x,y)\in\bR^d \times \bR^d : x=y\} \rightarrow \bR$ satisfying that there exists a positive constant $r_0$ such that
\begin{equation}
M_0 := \sup_{x\in\bR^d} \int_{\{y \in \R^d:|y-x| \ge r_0\}} J(x,y) dy <\infty, \label{e:intJ}
\end{equation}
and
\begin{equation}
 J(x,y) = \frac{f(|x-y|^{-1})}{|x-y|^d} \qquad\text{ for }\quad |x-y|<r_0,\label{e:J}
\end{equation}
where $f$ is a non-decreasing function from $[0,\, \infty)$ to $[0,\,\infty)$ having the following upper and lower growth  conditions at infinity;
\begin{equation}
a_1s^{\delta_1} \le \frac{f(st)}{f(t)} \le a_2s^{\delta_2} \qquad\text{ for }\quad s\ge 1,\, t\ge 1/r_0,\label{e:wsc}
\end{equation}
with some constants $a_1,\, a_2>0$ and $\delta_1,\, \delta_2\in(0,\, 2)$. 

If \eqref{e:J} holds for all $x,y \in \R^d$, then $J$ is the jump density of an isotropic unimodal L\'evy process. 
As an important example of Markov process having jump density comparable to $J$, isotropic unimodal L\'evy processes contain many processes, for example, subordinate Brownian motions which have been dealt with in a lot of literatures in probability and potential theory.
Recently, in \cite{BGR} it is shown that \eqref{e:wsc} holds if the characteristic exponent 
of the corresponding isotropic unimodal L\'evy process has the above upper and lower growth  conditions at infinity.

In this paper we fix  constants $\lambda,\, \Lambda>0$.
With these $\lambda,\, \Lambda>0$ and fixed function $J$ in \eqref{e:intJ} and \eqref{e:J}, 
 we define integro-differential operators comparable to $J$ as follows;
let kernel $K : \bR^d \times \bR^d\setminus \{(x,y)\in\bR^d \times \bR^d : x=y\} \rightarrow \bR$ be comparable to $J$ uniformly as 
\begin{equation} 
\lambda J(x,y) \le K(x,y) \le \Lambda J(x,y), \quad x,y \in \bR^{d}. \label{a:com}
\end{equation}
Define integro-differential operator with the kernel $K$ as 
\begin{equation}
L u(x) = L_K u(x) =\int_{\bR^d} \left(u(y)-u(x)-\nabla u(x)\cdot (y-x) \1_{\{|y-x|<r_0\}}\1_{\{\delta_2\ge 1\}}\right)K(x,y) dy. \label{d:lop} 
\end{equation}
Note that $K$ may not be symmetric. 

We will investigate the H\"older regularity of the solution to fully nonlinear nonlocal equation $\sI u = 0$ in an open set $D$ which is uniformly elliptic with respect to a family $\sL$ of integro-differential operators $L$, i.e.,
\begin{equation}\label{e:uniformly elliptic}
\sM^{-}_{\sL} (u-v) (x) \le \sI u(x) -\sI v(x) \le \sM^{+}_{\sL} (u-v) (x) \quad\text{ for } x\in  D,
\end{equation}
where maximal and minimal operators $\sM^{+}_{\sL}$, $\sM^{-}_{\sL}$ are defined by 
$$ \sM^{+}_{\sL} u(x) = \sup_{L\in\sL} Lu(x) \quad \text{ and } \quad \sM^{-}_{\sL} u(x) = \inf_{L\in\sL} Lu(x).$$
In \cite{CS09} Caffarelli and Silvetre established the theory which can be used to prove the Harnack inequality and the H\"older continuity of the solution $u$ with respect to the family of integro-differential operators whose kernels $K(x,y)$'s satisfying
\begin{equation} \label{e:CS}
 \frac{(2-\alpha)\lambda}{|x-y|^{d+\alpha}} \le K(x,y) \le \frac{(2-\alpha)\Lambda}{|x-y|^{d+\alpha}}, \quad x,y\in \bR^d,
\end{equation}
and
\begin{equation*}
K(x,x+z)=K(0,z)=K(0,-z), \quad x,z\in\bR^d. 
\end{equation*}

The purpose of this paper is to obtain the H\"older continuity of solutions to uniformly elliptic fully nonlinear nonlocal equation with respect to more general class of operators. Two classes of integro-differential operators in our consideration are
\begin{equation*}
 \sL := \left\{ L=L_{K} :  K \text{ satisfying } \eqref{a:com} \right\},
\end{equation*}
and
\begin{equation*}
 \sL_{sym} := \left\{ L_{K}\in \sL :  K(x,x+z)=K(x,x-z)\text{ for } x,z\in\bR^d \right\} .
 \end{equation*}
Uniformly ellipticity of the operator $\sI$ in \eqref{e:uniformly elliptic} allows us to regard the solution to $\sI u = 0$ as the sub and super solution to maximal and minimal operators respectively. 

The main result in this paper is following. See Section \ref{S:vis} for  the definition of solutions in the viscosity sense. 
\begin{thm}\label{thm1}
Assume that  \eqref{e:intJ}, \eqref{e:J} and \eqref{e:wsc} hold for the function $J(x,y)$ with $\delta_1 \in (1,2)$ or $\delta_2 \in (0,1)$. Let $z_0 \in \R^d$, $r>0$ and $u$ be a bounded function so that 
\begin{eqnarray*}
\sM^{+}_{\sL} u(x) &\ge& 0 \qquad for \quad x \in B(z_0,r) \\
\sM^{-}_{\sL} u(x) &\le& 0 \qquad for \quad x \in B(z_0,r)
\end{eqnarray*}
in the viscosity sense. Then there exist two constants $\alpha>0$ and $C>0$ such that
\begin{equation*}
\sup_{x,y \in B(z_0,{r/2})} \frac{|u(x)-u(y)|}{|x-y|^{\alpha}} \le C\left(\frac{1}{r\wedge r_0}\right)^{\alpha}\Arrowvert u \Arrowvert _{\infty}
\end{equation*}
where the values of $\alpha$ and $C$ depend only on $d,r_0,\lambda,\Lambda,a_1,a_2,\delta_1,\delta_2$ and $f(r_0^{-1})$.

In the case $\delta_1\le 1\le \delta_2$, the above assertion still holds for $\sL_{sym}$ in place of $\sL$.
\end{thm}

In \cite{S06} Silvestre developed an analytic method to obtain this result for integro-differential operator $L$ satisfying an assumption that given $\delta>0$ and some auxiliary function $b(x) = (1-|x|^2)_+^2$, there are positive constants $\kappa$ and $\eta$ such that for every $x\in\bR^d$
\begin{equation}\label{a:sil}
\kappa Lb(x) + 2\int_{|y-x|\ge 1/4} (|8(y-x)|^{\eta}-1)K(x,y) dy < \frac{1}{2} \inf_{A\subset B(0,2), |A|>\delta} \int_{A} K(x,y) dy.
\end{equation}
In \cite{S06} classes of operators with similar behavior at every scale were considered.
Thus, by considering an operator $L_{r,x_0}$ obtained by a change of variables,  Silvestre obtained the H\"older continuity of the solution $u$ to $Lu=0$ in $B(0,1)$ using the fact that \eqref{a:sil} holds for every $L_{r,x_0},\  r>0, x_0\in\bR^d$ by the scaling. 
Arguments using scaling assumption was employed to prove the same result for fully nonlinear equation with respect to the family of operators having kernels comparable to that of fractional Laplacian.

In this paper we will deal with the kernels comparable to $J(x,y)$ which is equal to $f(|x-y|^{-1})|x-y|^{-d}$ on near diagonal part. Our $J(x,y)$ has no scaling property and  it is not comparable to fractional Laplacian near zero.  Instead, we assume the weak growth conditions \eqref{e:wsc} of $f$ instead of the stability of the operator like fractional Laplacian. In this paper, we extend \eqref{a:sil} 
to our $J$ and for all small scales. 

Caffarelli and  Silvestre established the regularity results for fully nonlinear nonlocal equation extending the method for differential equation in \cite{CS09}. By this work they could obtain the result uniformly as $\alpha$ in \eqref{e:CS} goes to $2$. But they made the best use of the symmetry of the kernels which imply no effect of the gradient term in integro-differential operators. In our consideration we just impose the symmetry condition in the case $\delta_1\le1 \le \delta_2$. There are some results for the non symmetric case when the symmetric part of operators are comparable to fractional Laplacian (see \cite{KL12, KL13, CD12, KRS}).

Our paper is organized as follows. In Section \ref{S:vis} we recall the notion of viscosity solution in \cite{CS09}. 
 In Section \ref{S:pf_result} we prove the main result by modifying the method in \cite{S06} properly to our consideration. 
In Section \ref{S:iulp} we give examples covered in this paper. We recall recent results in \cite{BGR} on 
the estimates on densities of isotropic unimodal L\'evy processes, which serve as main examples.   

Throughout this paper, $d\geq 1$ and  the constants $r_0$, $M_0$, $\lambda$, $\Lambda$, $a_1$, $a_2$, $\delta_1$ and $\delta_2$ will be fixed. We use $C_1, C_2, C_3$ to denote the constants which are obtained in the proofs of theorems and depends only the aforementioned constants. We define $\osc_E u = \sup_{x\in E} u(x) - \inf_{x\in E} u(x)$ for a subset $E$ in $\bR^d$. 
For any Borel subset $E\subset\bR^d$, 
$\overline{E}$ stands for the closure of $E$. We denote by $\omega_d= 2\pi^{d/2}/\Gamma(d/2)$ the surface measure of the unit sphere in $\bR^d$.

\section{Preliminaries: Viscosity solution and Elliptic operator}\label{S:vis}
The notion of viscosity solution is very useful to solve Bellman, Hamilton-Jacobi-Bellman or Isaacs equations. This notion allows us to apply smooth functions to an operator instead of measurable functions whenever there is a test function touching from above or below. For integro-differential operators  we follow the definition of viscosity solution described in \cite{CS09}. 

\begin{defn}
A function $u : \bR^d \rightarrow \bR$, upper(resp. lower) semicontinuous in $\overline{D}$, is said to be a subsolution(resp. supersolution) to $\sI u = g$, and we write $\sI u \ge g$(resp. $\sI u \le g$) if the following holds : If we have a function $v$ defined by 
\begin{equation*}
v(y) = \begin{cases}
\varphi(y) & \text{ if } y \in N_x \\
u(y) & \text{ if } y \in \bR^d\setminus N_x
\end{cases}
\end{equation*}
where $N_x$ is a neighborhood of $x$ in $D$ and $\varphi$ is a $C^2$ function touching $u$ from above (resp. below) at $x$, i.e. $\varphi (x) = u(x)$ and $\varphi(y) >u(y)$ (resp. $\varphi(y) < u(y)$) for $y \in N_x\setminus\{x\}$,
then $\sI v(x) \ge g(x)$ (resp.  $\sI v(x) \le g(x)$).
A solution is a function $u$ that is both a subsolution and a supersolution.
\end{defn}

Since Jensen developed the idea to find uniqueness results of viscosity solutions to fully nonlinear second order differential equations which is of the form $F(D^2u,Du,u)=0$ in \cite{J88}, Ishii established in \cite{Ish89} a lemma, so called Jensen-Ishii's lemma, which is used to obtain comparison principles for fully nonlinear degenerate elliptic equations. He obtained the existence and uniqueness of the continuous solution by combining the comparison principle and Perron's method.

 The existence and uniqueness of the solution for the nonlocal Dirichlet problem 
 \begin{align*}
 &\sI u(x) = 0 \quad \text{ in } D \\
 &u(x) = g(x) \quad \text{ for } x \in \bR^d\setminus D
\end{align*}
have been studied in a lot of literatures. In \cite{A91a, A91b} Awatif successfully applied the Jensen's method to obtain the results for first order differential equations with an integro-differential term. Barles and Imbert gave a general proof for second order elliptic integro-differential equations in which a new definition of viscosity solutions equivalent to the  above definition is introduced (see \cite[Definition 4]{BI08}). In \cite{CS09} Caffarelli and Silvestre considered a somewhat abstract class of operators, i.e., elliptic operators (see Definition \ref{d:ell}), and established the general method used to find the unique viscosity solution for the translation invariant uniformly elliptic nonlocal equations. They obtained the comparison principle for elliptic operator $\sI$ of supremum or infimum type. In \cite{BCI} Barles, Chasseigne and Imbert developed the viscosity solution theory for nonlocal nonlinear equations in quite different assumptions from \cite{CS09}. They assumed that the kernels have a certain continuity in $x$. Although the operators in this paper are motivated from these studies, the existence and uniqueness results has not been found 
yet because it has considerably general kernels in the integro-differential operators.

We define a class of functions $C^{1,1}(x)$ as follows: a function $u$ is said to be in $C^{1,1}(x)$ if there is a vector $v\in\bR^d$ and constants $M,\epsilon>0$ such that $ |u(y)-u(x)-v\cdot(y-x)|\le M|y-x|^2 $ for $|y-x|<\epsilon$. 

The following is the minimal condition to obtain the comparison principle for various equations.

\begin{defn}\label{d:ell}
Let $\sL$ be a class of linear integro-differential operators.  An elliptic operator $\sI$ with respect to $\sL$ is an operator with the following properties:
\begin{itemize}
\item[$\bullet$] If $u$ is any bounded function, $\sI u(x)$ is well-defined for all $u\in C^{1,1}(x)$.
\item[$\bullet$] If $u$ is $C^2$ in some open set $D$, then $\sI u(x)$ is a continuous function in $D$.
\item[$\bullet$] If $u$ and $v$ are bounded functions $C^{1,1}$ at $x$, then
$$
\sM^{-}_{\sL} (u-v) (x) \le \sI u(x) - \sI v(x) \le \sM^{+}_{\sL} (u-v)(x).
$$ 
\end{itemize}
\end{defn}

The last condition in the above definition is used to linearize the equation through the extremal operators. A study on fully nonlinear elliptic differential equations is in \cite{CC} and references therein. For the fully nonlinear uniformly elliptic \emph{nonlocal} operator Caffarelli and Silvestre systematically established the Harnack inequality, the H\"older estimate and $C^{1,\alpha}$ regularity of the solutions to the Dirichlet problems in \cite{CS09}. They dealt with translation invariant elliptic operator with respect to the symmetric kernels comparable to the jumping kernel of an isotropic $\alpha$-stable process. They developed the ABP-estimate for integro-differential operators. After that, in \cite{CS11} the same authors obtained the regularity results for nonlocal elliptic equations with kernels which are not translation invariant. They used the closeness of the operator with translation invariant one which already has regularity results. After their works, operators having non-symmetric kernels were dealt with in \cite{CD12, KRS, KL12, KL13} . In this case a difficulty comes from an effect of the gradient.

\section{The proof of Theorem \ref{thm1}} \label{S:pf_result}
We mainly follow the method developed by Silvestre in \cite{S06}. 
Recall that 
 $a_1, a_2, \delta_1, \delta_2, r_0, \lambda, \Lambda, M_0$ 
 are fixed constants in \eqref{e:intJ}, \eqref{e:J}, \eqref{e:wsc} and \eqref{a:com}. 

 We start from simple calculations for integrals of a non-decreasing function satisfying local growth conditions.

\begin{lemma}
Let $f$ be a non-decreasing function from $[0,\infty)$ to $[0,\infty)$ satisfying \eqref{e:wsc}. Then we have the following, for $0<r<r_0$,
\begin{align}
&r^{-2} \int_0^r s f(s^{-1}) ds \le \frac{a_2}{2-\delta_2}f(r^{-1}),\label{e:f0}\\
&\int_r^{r_0} s^{-1} f(s^{-1}) ds \le \frac{1}{a_1\delta_1}f(r^{-1}),\label{e:fr0}\\
& r^{-1} \int_0^rf(s^{-1}) ds \le \frac{a_2}{1-\delta_2}f(r^{-1}), & \text{ if } \delta_2<1,\label{e:gradf0}\\
& r^{-1} \int_r^{r_0}f(s^{-1}) ds \le \frac{1}{a_1(\delta_1-1)} f(r^{-1}), &\text{ if } \delta_1>1.\label{e:gradfr0}
\end{align}
\end{lemma}

\pf Since we have $s<r<r_0$, we get $f(s^{-1})\le a_2r^{\delta_2}s^{-\delta_2}f(r^{-1})$ by the upper growth condition in \eqref{e:wsc}. Therefore we obtain
\begin{equation*}
r^{-2} \int_0^r  s f(s^{-1}) ds \le a_2 r^{\delta_2-2} f(r^{-1}) \int_0^r s^{1-\delta_2} ds = \frac{a_2}{2-\delta_2}f(r^{-1}),
\end{equation*}
and
\begin{equation*}
r^{-1}\int_0^r  f(s^{-1}) ds \le a_2 r^{\delta_2-1} f(r^{-1}) \int_0^r s^{-\delta_2} ds = \frac{a_2}{1-\delta_2}f(r^{-1}), \qquad \text{ for } \delta_2<1.
\end{equation*}
We have proved  \eqref{e:f0} and \eqref{e:gradf0}. Using the lower growth condition in \eqref{e:wsc}, the proofs of inequalities \eqref{e:fr0}  and \eqref{e:gradfr0} are similar and we omit the proof. \qed

\begin{lemma}\label{l:wedge}
Suppose $J(x,y)$ satisfies \eqref{e:intJ}, \eqref{e:J} and \eqref{e:wsc}. Then we have 
\begin{equation}
\sup_{x \in \R^d}\int_{\bR^d} \left( 1\wedge\left(\frac{|y-x|}{r}\right)^2 \right) J(x,y) \,dy \le  C_1 f(r^{-1}) \quad \text{ for }\quad  r \le r_0, \label{e:wedge}
\end{equation}
where $C_1 =\omega_d \left( \frac{a_2}{2-\delta_2} + \frac{1}{a_1\delta_1} \right) + \frac{M_0}{f(r_0^{-1})}.$
\end{lemma}
\pf We first decompose the integral in the left side of \eqref{e:wedge} into three parts and 
 use \eqref{e:intJ} and \eqref{e:J} so that
\begin{align*} 
&\int_{\bR^d} \left( 1\wedge\frac{|y-x|^2}{r^2} \right) J(x,y)\, dy\\
=&\int_{|y-x|<r} \frac{|y-x|^2}{r^2} \frac{f(|y-x|^{-1})}{|y-x|^d}\, dy + \int_{r\le |y-x|<r_0} \frac{f(|y-x|^{-1})}{|y-x|^d} \,dy + \int_{|y-x|\ge r_0} J(x,y)\, dy \\
\le&\  \omega_d \int_0^r r^{-2}sf(s^{-1}) \, ds + \omega_d \int_r^{r_0} s^{-1}f(s^{-1}) \, ds + M_0.
\end{align*}
By \eqref{e:f0}, \eqref{e:fr0} and the monotonicity of $f$, we have \eqref{e:wedge}.
\qed

When the integro-differential operator $L$ is the fractional Laplacian, the effect on $Lu$ of scaling to enlarge $u$ is transferred to the the scaling of kernel. In the next lemma we show that,  for our operator $L$,  
the effect of  magnifying support outside of 
the origin can be controlled by the growth conditions of the kernel near the origin. 
\begin{lemma}\label{l:u_growth}
For any $\epsilon >0 $ there are constants $r=r(\epsilon)\in (0,r_0)$ and $\eta=\eta(\epsilon)\in(0,\delta_1)$ such that for all $0<s<r$ 
\begin{equation}\label{e:mag_sup_out}
\sup_{ x \in \bR^d} \int_{|y-x|>\frac{s}{4}} \left(\left(2\frac{|4(y-x)|\wedge r_0}{s}\right)^{\eta}-1\right) J(x,y) dy < \epsilon f(s^{-1})
\end{equation}
\end{lemma}
\pf
We decompose the integral in the left side of \eqref{e:mag_sup_out} into two parts according to the distance between $x$ and $y$ as following
\begin{align*}
&\quad\int_{|y-x|>s/4} \left(\left(2\frac{|4(y-x)|\wedge r_0}{s}\right)^{\eta}-1\right) J(x,y) \, dy \\
& =\int_{|y-x|>r_0/4} ( 2^{\eta}(r_0/s)^{\eta} - 1) J(x,y) \, dy   + \int_{s/4<|y-x|\le r_0/4} \left(2^{\eta}\left(\frac{4|y-x|}{s}\right)^{\eta}-1\right) J(x,y) \, dy \\
 & =: I_1 + I_2.
\end{align*}
From the obvious inequality 
$
\1_{\{|x|\ge r_0/4\}} \le 1\wedge (4|x|/r_0 )^2\le 16 (1\wedge (|x|/r_0))^2,
$ and  \eqref{e:wedge}
we obtain the bound for $I_1$ as  
\begin{align*}
I_1 \le 16(2^{\eta}(r_0/s)^{\eta}-1)\int_{\bR^d} \left(1\wedge\left(\frac{|y-x|}{r_0}\right)^2\right) J(x,y) dy \le 2^{4+\eta}(r_0/s)^{\eta}C_1f(r_0^{-1}).
\end{align*}
Since, by the lower growth condition of \eqref{e:wsc}, $f(s^{-1})\ge a_1 (r_0/s)^{\delta_1}f(r_0^{-1})$, if $\eta$ is less than $\delta_1$ then 
$$
I_1 \le 2^{4+\delta_1}a_1^{-1}C_1(s/r_0)^{\delta_1-\eta}f(s^{-1}).
$$

On the other hand, by \eqref{e:J}
\begin{align*}
I_2 & =  \int_{s/4<|y-x|\le r_0/4}  \left(2^{\eta}\left(\frac{4|y-x|}{s}\right)^{\eta}-1\right) \frac{f(|y-x|^{-1})}{|y-x|^d} dy \\
& = \omega_d \int_s^{r_0} \left(2^{\eta}(t/s)^{\eta}-1\right)\frac{f(4/t)}{t} dt
\end{align*}
Since $f(4t^{-1}) \le a_2 4^{\delta_2} f(t^{-1}) \le a_2 4^{\delta_2}a_1^{-1} (t/s)^{-\delta_1}f(s^{-1})$ by \eqref{e:wsc}, 
we have 
\begin{align*}
I_2 &\le \omega_d  a_2 4^{\delta_2} a_1^{-1} f(s^{-1}) \int_s^{r_0} (2^{\eta}(t/s)^{\eta} -1)(t/s)^{-\delta_1}t^{-1} dt  \\
&=\omega_d  a_2 4^{\delta_2} a_1^{-1} f(s^{-1}) \int_{1}^{r_0/s}( (2t)^{\eta} -1) t^{-\delta_1-1} dt \\
&\le \omega_d a_2 4^{\delta_2} a_1^{-1}f(s^{-1})\int_{1}^{\infty}( (2t)^{\eta} -1) t^{-\delta_1-1} dt 
\end{align*}
By the dominated convergence theorem, we can choose $\eta<\delta_1$ such that $\omega_d a_2 4^{\delta_2} a_1^{-1} \int_1^{\infty} ( (2t)^{\eta} -1) t^{-\delta_1-1} dt<\eps/2$ and then find $r$ such that $ 2^{4+\delta_1}a_1^{-1}C_1(r/r_0)^{\delta_1-\eta} < \eps/2$.
\qed

 Define test functions 
 $$
\beta(t) = (1-t^2)_+^2, \quad t\ge 0, \quad \text{ and } \quad  b_{z,r}(x) := \beta(|x-z|/r), \quad x,z\in\bR^d, r>0.
$$
 In the following lemma we prove that the minimal operator applied to test function $b_{z,r}$  is bounded by $c f(r^{-1})$ for 
 all $r<r_0$. In \cite{S06} Silvestre obtained this result for the operators like fractional Laplacian for every $r>0$ using scaling.

\begin{lemma}\label{l:u_bump}
For any $0<r\le r_0$ and $x, z \in \bR^d$, if $\delta_1\in(1,2)$ or $\delta_2\in(0,1)$ then we have  
\begin{equation*}
\left| \sM^{-}_{\sL}b_{z,r}(x) \right| \le C_2 f(r^{-1}) 
\end{equation*}
where the constant $C_2$ is depending only on $d, a_1,a_2,\delta_1,\delta_2, \Lambda, M_0$ and $f(r_0^{-1})$. In the case $\delta_1\le 1\le \delta_2$ we have same bound for $\sM^{-}_{\sL_{sym}}b_{z,r}(x)$ with $12 d \Lambda C_1$ in place of $C_2$.
\end{lemma}
\pf We consider the following three integrals;
\begin{align*}
I_1 &= \int_{\bR^d} \left( b_{z,r}(y)-b_{z,r}(x)-\nabla b_{z,r}(x)\cdot (y-x)\1_{\{|y-x|<r\}} \right) K(x,y) dy , \\
I_2 &= \int_{ |y-x|<r} \left( \nabla b_{z,r}(x)\cdot (y-x) \right) K(x,y) dy , \\
I_3 &= \int_{r\le |y-x|<r_0} \left(\nabla b_{z,r}(x)\cdot (y-x) \right) K(x,y) dy.
\end{align*}
Since $Lb_{z,r}(x) =L_Kb_{z,r}(x) = I_1+I_2\1_{\{\delta_2<1\}}-I_3\1_{\{\delta_2\ge 1\}}$,  it is enough to estimate $I_1$, $I_2$ and $I_3$.

First, from the definition of $b_{z,r}$ we have that for $|y-x|<r$
\begin{align*}
&
| b_{z,r}(y)-b_{z,r}(x)-\nabla b_{z,r}(x)\cdot (y-x)\1_{\{|y-x|<r\}} |\\
\le & \sup_{w\in\bR^d, 1\le i,j\le d} \left|\frac{\partial^2b_{z,r}}{\partial x_i \partial x_j}(w)\right| d|y-x|^2 \le \frac{d|y-x|^2}{r^2}\sup_{0\le s \le 1} \Big\{|\beta''(s)|+|\beta'(s)/s|\Big\} \le 12d\frac{|y-x|^2}{r^2},
\end{align*}
and clearly the integrand in $I_1$ is  bounded by $K(x,y)$ for $|y-x|\ge r$. Thus by \eqref{a:com}
and  \eqref{e:wedge} 
\begin{align}\label{e:I_1}
|I_1|\le 12d\Lambda \int_{\bR^d} \left( 1\wedge\frac{|y-x|^2}{r^2} \right) J(x,y) \, dy \le 12d\Lambda C_1 f(r^{-1}).\end{align}

The gradient term in the integrand of $I_2$ and $I_3$ is bounded by $\sup_{0\le s\le 1} |\beta'(s)|\le 4$. Thus, by \eqref{e:J} and \eqref{e:gradf0}, we have the bound for $I_2$  as
\begin{align}\label{e:I_2}
|I_2| &\le 4\Lambda \int_{|y-x|<r} \frac{|y-x|}{r} J(x,y)\, dy = 4\Lambda \int_{|y-x|<r} \frac{|y-x|}{r} \frac{f(|y-x|^{-1})}{|y-x|^d}\, dy \nn\\
&= 4\omega_d\Lambda r^{-1}\int_0^r f(s^{-1}) \, ds \le 4\omega_d\Lambda \frac{a_2}{1-\delta_2}f(r^{-1}), \quad \text{ if } \delta_2<1.
\end{align}
By using \eqref{e:gradfr0} instead of \eqref{e:gradf0} we have the bound for $I_3$  as
\begin{align}\label{e:I_3}
|I_3| &\le 4\Lambda \int_{r\le|y-x|<r_0} \frac{|y-x|}{r} J(x,y) dy = 4\Lambda  \int_{r\le|y-x|<r_0} \frac{|y-x|}{r} \frac{f(|y-x|^{-1})}{|y-x|^d}  dy\nn \\
 &= 4\omega_d\Lambda \int_r^{r_0}r^{-1}f(s^{-1}) ds \le 4\omega_d\Lambda  \frac{1}{a_1(\delta_1-1)}f(r^{-1}),\quad \text{ if } \delta_1>1.
\end{align}

Combining \eqref{e:I_1}--\eqref{e:I_3} and then taking $\inf_{L\in\sL}$, we get 
that for $r\le r_0$ if $\delta_2<1$ or $\delta_1>1$, 
$\left| \sM^{-}_{\sL}b_{z,r}(x) \right| \le C_2 f(r^{-1})$  
where $C_2 = 12d\Lambda\left(C_1 + \omega_d \frac{a_2}{1-\delta_2} + \omega_d \frac{1}{a_1(\delta_1-1)}\right)$.

In the case $\delta_1\le 1 \le \delta_2$ we can get $Lb_{z,r}(x) = I_1$ for $L \in \sL_{sym}$ because $I_3=0$ by the symmetry of $K(x,x+\cdot)$. Therefore we obtain the bound $|\sM^{-}_{\sL_{sym}}b_{z,r}| \le 12 d \Lambda C_1 f(r^{-1})$.
\qed

In the proof of the following theorem we will consider a test function $\psi$, in 
the definition of viscosity solution,  touching $u$ from above at a maximum point
. 
\begin{thm}\label{l:main}
Suppose $\delta_1\in(1,\, 2)$ or $\delta_2\in(0,\, 1)$. Then there exist constants $r_1\in(0,r_0)$ and $\eta_1>0$ such that if $u$ is a function that satisfies the following assumptions for $z\in \bR^d$, $0<r<r_1$
\begin{alignat}{2}
\sM^{+}_{\sL}u(x) &\ge 0 && \qquad \textrm{for } x \in B(z,r), \nn\\
u(x) & \le \frac12 && \qquad \text{for } x \in B(z,r), \nn\\ \label{e:growth}
u(x) & \le \left(2\frac{|x-z|\wedge r_0}{r}\right)^{\eta_1} -\frac{1}{2} && \qquad \text{for } x \in \bR^d\setminus B(z,r), \\
\frac12|B(z,r)|&< | \{ x \in B(z,r) : u(x) \le 0\} |, &&  \nn
\end{alignat}
then $u \le 1/2 -\gamma$ in $B(z,r/2)$ for some constant $\gamma \in (0, 1-2^{-\eta_1})$ depending on $r_1$ and $\eta_1$.

If we suppose $\delta_1\le 1\le \delta_2$ then the above assertion holds for $\sL_{sym}$ instead of $\sL$.
\end{thm}

\pf 
The proofs are the same for the case $\delta_1\in(1,\, 2)$ or $\delta_2\in(0,\, 1)$ and  the case $\delta_1 \le 1 \le \delta_2$. So we give the proof for the case $\delta_1>1$ or $\delta_2<1$ only.

We first choose $r_1\in(0,r_0/2)$ and $\eta_1>0$: By Lemma \ref{l:u_growth} with $\epsilon=\frac{ \omega_d\lambda}{\Lambda 2^{d+5+\delta_2}a_2d}$ there exist
$r_1\in(0,r_0)$ and $\eta_1>0$
such that for any $r<r_1$
\begin{equation}\label{ee:u_growth}
\int_{|y-x_1|\ge\frac{r}{4}} \left(\left(2\frac{|4(y-x_1)|\wedge r_0}{r}\right)^{\eta_1} -1 \right) J(x_1, y) \, dy  \le \frac{ \omega_d\lambda}{\Lambda 2^{d+5+\delta_2}a_2d}f(r^{-1})
\end{equation}

$\theta$ is  a small positive constant depending on $r_1$ and $\eta_1$ which will be chosen later. Define $\gamma = \theta(\beta(1/2)-\beta(3/4))$.

Suppose there is a point $x_0 \in B(z,r/2)$ such that $u(x_0) > 1/2-\gamma = 1/2-\theta\beta(1/2)+\theta\beta(3/4)$. Then we have
\begin{equation*}
u(x_0)+\theta b_{z,r}(x_0)\ge u(x_0)+\theta \beta(1/2) > 1/2 + \theta \beta(3/4) \ge u(x) + \theta b_{z,r}(x) 
\end{equation*}
for  $ x\in B(z,r)\setminus B(z,3r/4)$. This means that the supremum of $u+\theta b_{z,r}$ in $B(z,r)$ is greater than $1/2$ and is taken at an interior point $x_1$ of $B(z,3r/4)$. Since $u+\theta b_{z,r}$ has a maximum at $x_1$, we have a test function $\varphi$ touching $u+\theta b_{z,r}$ from above at $x_1$. For $\tilde\epsilon >0$ and $0<s<r/4$ define a function $\varphi$ by 
\begin{displaymath}
\varphi(y) = 
\begin{cases}
u(x_1)+\theta b_{z,r}(x_1)+\tilde\epsilon |y-x_1|^2 & \text{ if } |y-x_1|<s, \\
u(y)+\theta b_{z,r}(y) & \text{ otherwise. } 
\end{cases}
\end{displaymath}

Now we evaluate $\sM^{+}_{\sL}\varphi(x_1)$. On the one hand, the fact that $\varphi - \theta b_{z,r}$ is a test function touching $u$ from above at $x_1$ and the ellipticity of $\sM^{+}_{\sL}$ imply 
\begin{equation}\label{e:var_b}
\sM^{+}_{\sL}\varphi(x_1) \ge \sM^{+}_{\sL}(\varphi-\theta b_{z,r})(x_1) + \theta \sM^{-}_{\sL}b_{z,r}(x_1) \ge \theta \sM^{-}_{\sL}b_{z,r}(x_1).
\end{equation}
On the other hand, from $\nabla \varphi(x_1) =0$  we can divide $L\varphi(x_1)$ into two parts as follows
\begin{align}
L\varphi(x_1)& =  \int_{|y-z|\ge r} \left( \varphi(y) - \varphi(x_1) \right)K(x_1,y)\, dy 
 + \int_{|y-z|< r} \left( \varphi(y) - \varphi(x_1) \right)K(x_1,y)\, dy \nn \\
& =: I_1+I_2. \label{e:n1}
\end{align}
Since the support of $b_{z,r}$ is in $B(z,r)$ and $\varphi(x_1)>1/2$, by (\ref{e:growth}) and (\ref{a:com})
\begin{equation*}
I_1 =  \int_{|y-z|\ge r } (u(y)-\varphi(x_1) ) K(x_1,y) \,dy \le \int_{|y-z|\ge r} \left(\left(2\frac{|y-z|\wedge r_0}{r}\right)^{\eta_1} - \frac{1}{2}-\frac{1}{2}\right)\Lambda J(x_1,y) \,dy. 
\end{equation*}
Since  $|y-x_1|\ge |y-z|-|z-x_1|> r/4$ and $|x_1-z|\le 3|y-x_1|$ if $|x_1-z|<3r/4$ and $|y-z|\ge r$,  we have
\begin{equation}
I_1 \le \Lambda\int_{|y-x_1|\ge\frac{r}{4}} \left(\left(2\frac{|4(y-x_1)|\wedge r_0}{r}\right)^{\eta_1} -1 \right) J(x_1,y) \, dy.\label{e:n2}
\end{equation}
To estimate $I_2$ we decompose the region $\{y:|y-z|<r\}$ in the integral  into $\{y: |y-z|<r ,\, |y-x_1|<s \}$ and $\{ y: |y-z|<r ,\, |y-x_1|\ge s \}$ and control the integrand as 
\begin{align}
I_2 & = \int_{\substack{|y-z|<r \\ |y-x_1|<s}}(\varphi(y)-\varphi(x_1))K(x_1,y) dy + \int_{\substack{|y-z|<r \\ |y-x_1|\ge s}} (\varphi(y)-\varphi(x_1))K(x_1,y) dy \nn\\
&\le \int_{|y-x_1|<s} \tilde\epsilon |y-x_1|^2 \Lambda J(x_1,y) dy +\int_{\substack{|y-z|<r \\ |y-x_1|\ge s}} (u(y)+\theta b_{z,r}(y)-u(x_1)-\theta b_{z,r}(x_1))K(x_1,y) dy. \nn
\end{align}
The first term is bounded by $\tilde\epsilon\Lambda r_0^2 C_1f(r_0^{-1})$ for $s<r_0$ by \eqref{e:wedge}. The integrand in the second term is nonpositive so that the value of the integral will be greater if we restrict the region in the integral to the set $A_s := \{y\in B(z,r) : |y-x_1|\ge s \text{ and } u(y) \le 0 \}$. If $s$ approaches $0$ then $|A_s|$ goes to $|\{y\in B(z,r) : u(y) \le 0 \}|$ which is greater than $|B(z,r)|/2$. So we can find small $s>0$ satisfying $A_s \subset B(x_1,2r)$ and $|A_s| > |B(z,r)|/2$. Thus we have the following estimate for $\theta < 1/4$
\begin{align}
I_2 &\le \tilde\epsilon\Lambda r_0^2 C_1f(r_0^{-1})+ (\theta - \frac{1}{2})\lambda  \int_{\substack{|y-z|<r, |y-x_1|\ge s \\ u(y)\le 0 }}J(x_1,y) dy \nn \\
&\le  \tilde\epsilon\Lambda r_0^2 C_1f(r^{-1}) -\frac{\lambda}{4} \inf_{\substack{A\subset B(x_1,2r),\\ |A| >|B(z,r)|/2}} \int_A  J(x_1,y) dy\,. \nn
\end{align}
For the last term above we observe that, by \eqref{e:J}, the monotonicity of $f$ and \eqref{e:wsc},  
we have 
for $r<r_0/2$ and $A \subset B(x_1,2r)$ with $|A|>|B(z,r)|/2$, 
\begin{align*}
\frac{\lambda}{4}\int_A  J(x_1,y) dy &= \frac{ \lambda}{4}\int_A \frac{f(|y-x_1|^{-1})}{|y-x_1|^d} dy \ge \frac{ \lambda}{4}\int_A  \frac{f((2r)^{-1})}{2^d r^d} dy\\
	&\ge  \frac{ \lambda}{2^{d+2}}a_2^{-1}2^{-\delta_2}f(r^{-1})r^{-d}|A|  \ge  C_3 f(r^{-1})
\end{align*}
where $\displaystyle C_3 = \frac{ \omega_d\lambda}{2^{d+3+\delta_2}a_2d}$.
Thus we have $ I_2\le -C_3 f(r^{-1})/2$ for $\theta<1/4$ and $\tilde\epsilon <C_3/(2\Lambda r_0^2 C_1)$. 
Hence from 
\eqref{e:n1}, \eqref{e:n2} and this, 
we obtain the following 
\begin{equation*}
\sM^{+}_{\sL}\varphi(x_1) \le \Lambda\int_{|y-x_1|\ge\frac{r}{4}} \left(\left(2\frac{|4(y-x_1)|\wedge r_0}{r}\right)^{\eta_1} -1 \right) J(x_1, y) \, dy - \frac{C_3}{2}f(r^{-1})
\end{equation*}
for $\theta<1/4$ and $r<r_0/2$. Therefore, by Lemma \ref{l:u_bump}, \eqref{e:var_b} and \eqref{ee:u_growth}  for any $r<r_1$
\begin{equation}
-\theta C_2f(r^{-1}) \le \theta \sM^{-}_{\sL}b_{z,r}(x_1) \le \sM^{+}_{\sL}\varphi(x_1) \le -\frac{C_3}{4}f(r^{-1}). \nn
\end{equation}
We now choose $\theta =\frac{C_3}{8C_2}\wedge \frac{1-2^{-\eta_1}}{2(\beta(1/2)-\beta(3/4))}$ so that it yields a contradiction. \qed
\begin{remark}
We can prove the above \emph{Theorem \ref{l:main}} for $r^{d}\delta$ for some constant $\delta>0$ instead of $|B(z,r)|/2$.
\end{remark}

 When one deals with fraction Laplacian $\Delta^{\alpha/2}$, for example, the equation $\Delta^{\alpha/2}u=0$ in a ball $B(0,r)$, 
one may assume $r=1$ in the equation $Lu=0$ in a ball $B(0,r)$ by the scaling invariant property of the equation, i.e. $\Delta^{\alpha/2}\tilde{u} = 0$ in a ball $B(0,1)$ where $\tilde{u}(x) = u(rx)$.
 But in our case, we don't have such scaling property so we prove Theorem \ref{thm1} directly without using any scaling.

\vspace{0.5cm}

{\bf Proof of Theorem \ref{thm1}}.
Without loss of generality we assume that $z_0=0$. 
By normalization we can assume $\sup_{x\in\bR^d} |u(x)| = 1/2$. We have $r_1\in(0,r_0)$, $\eta_1>0$ and $\gamma \in (0, 1-2^{-\eta_1})$ such that if $u$ satisfy the assumptions in Theorem \ref{l:main} then $u>1/2-\gamma$. Let 
$\alpha := -\log_2(1-\gamma)$, which is less than $\eta_1$ and 
$x_0$ be a point in $B(0,r/2)$ and $s$ be the minimum of $r/2$ and $r_1/2$. 

We will show  by induction that 
\begin{align}\label{e:tclaim}
\osc_{B(x_0,2^{-k}s)}u \le (1-\gamma)^k \qquad \text { for all  } k = 0, 1, 2, \cdots. 
\end{align} 
 First, the case $k=0$ is true obviously. 
 
 Suppose $\osc_{B(x_0,2^{-k}s)}u \le (1-\gamma)^k$ for some nonnegative integer $k \ge 0$. Define $v(x) = (1-\gamma)^{-k}(u(x)-a_k)$ where $a_k = \min_{B(x_0,2^{-k}s)}u + (1-\gamma)^{k}/2$. We have two cases
\begin{align*}
\text{ (i) } & |\{x\in B(x_0,2^{-k}s) : v(x) \le 0\}|\ge |B(x_0,2^{-k}s)|/2; \\
\text{ (ii) } & |\{x\in B(x_0,2^{-k}s) : v(x) \ge 0\}|\ge |B(x_0,2^{-k}s)|/2.
\end{align*}
Without loss of generality  we assume (i) holds since  we may apply the same argument on $-v-(1/2-\max_{B(x_0,2^{-k}s)}v)$ for the  case (ii). 

Let's check the others conditions in Theorem \ref{l:main}. 
First, it is clear that 
$$\sM^{+}_{\sL}v(x) = (1-\gamma)^{-k}\sM^{+}_{\sL}u(x) \ge 0 \quad  \text{for } x \in B(x_0,2^{-k}s)$$ 
and 
$$v(x) \le (1-\gamma)^{-k}(\max_{B(x_0,2^{-k}s)}u - a_k) \le (1-\gamma)^{-k}(\osc_{B(x_0,2^{-k}s)}u - (1-\gamma)^{k}/2) \le 1/2.$$

We now check the third condition in Theorem \ref{l:main}. When $2^{-k+j}s\le |x-x_0| < 2^{-k+j+1}s$, $j = 0, \cdots, k-1$, we have 
\begin{align*}
v(x) &\le (1-\gamma)^{-k}\left(\max_{B(x_0,2^{-k+j+1}s)}u + \left(-\min_{B(x_0,2^{-k+j+1}s)}u + \min_{B(x_0,2^{-k}s)}u\right) -a_k\right) \\
 &\le (1-\gamma)^{-k}( (1-\gamma)^{k-j-1} -(1-\gamma)^k/2) \\
 &\le (1-\gamma)^{-\log_2({2|x-x_0|/(2^{-k}s)})} -\frac{1}{2}  \\
 &\le \left(\frac{2|x-x_0|}{2^{-k}s}\right)^{\eta_1} -\frac{1}{2}. 
\end{align*}
When $|x-x_0|\ge s$, we simply have 
\begin{displaymath}
v(x) \le (1-\gamma)^{-k}-\frac{1}{2} \le \left(2\frac{r_0}{2^{-k}s}\right)^{\eta_1} - \frac{1}{2}.
\end{displaymath}
Thus 
\begin{displaymath}
v(x) \le \left(2\frac{|x-x_0|\wedge r_0}{2^{-k}s}\right)^{\eta_1} - \frac{1}{2} \qquad\text{ for } |x-x_0| \ge 2^{-k}s.
\end{displaymath}
We have checked that all conditions in  Theorem \ref{l:main} holds. 
Therefore we obtain $v(x) \le 1/2 - \gamma$ for $|x-x_0| \le 2^{-k-1}s$, that is,  $\osc_{B(x_0,2^{-k-1}s)}u = (1-\gamma)^k\osc_{B(x_0,2^{-k-1}s)}v \le (1-\gamma)^{k+1}$. We have proved 
\eqref{e:tclaim}, which implies that $$|u(x)-u(x_0)| \le \left(2/s\right)^{\alpha}|x-x_0|^{\alpha} \le C (r\wedge r_0)^{-\alpha}|x-x_0|^{\alpha} \quad \text{ for all } x\in \bR^d$$ where $C= \left(4r_0/r_1\right)^{\alpha}$.
\qed

\section{Example: Isotropic unimodal L\'evy process}\label{S:iulp}
Let $(X_t, t\ge0)$ be a pure-jump isotropic L\'evy process in $\bR^d$.
Its characteristic function is
$$
\bE[\exp(i\xi\cdot  X_t)] = e^{-t\psi(|\xi|)}
$$
where  $\xi \to \psi(|\xi|)$ is called the characteristic exponent of $X$ and it has the representation 
$$
\psi(|\xi|) = \int_{\bR^d} (1-\cos(\xi\cdot x)) \nu(dx).
$$
The measure $\nu$ is  called the L\'evy measure of $X$ and it satisfies $\int(1 \wedge |x|^2) \nu(dx)< \infty$.
$(X_t, t\ge0)$ is called isotropic unimodal L\'evy process if the transition probability $\bP(X_t\in dx)$ has non-increasing density $p_t(x)$ with respect to Lebesgue measure. It is well known that $(X_t, t\ge0)$ is an isotropic unimodal L\'evy process if and only if 
the L\'evy measure $\nu(dx)$ of $X$ has non-increasing density, say, $\nu(x)$ (see \cite{W}). Note that $p_t(x)/t$ converges vaguely to $\nu(x)$. Denote $ \sup_{s\le t} \psi(s) $ by $\psi^ *(t)$. 

By \cite[Proposition 2]{BGR}  and \cite[Proposition 1]{G}, $\psi$ is almost increasing; 
\begin{equation}\label{e:almostinc}
\psi^* (t) \le \pi^2\psi(t) \quad \text{ for all }  t>0.
\end{equation}

Following upper bound  holds for $\nu$ without any extra condition (see \cite[Corollary 6]{BGR} and \cite[Theorem 2.2]{KKK}).
\begin{thm}For an isotropic unimodal L\'evy process $X$ in $\bR^d$, there is $C=C(d)$ such that 
\begin{equation}\label{nu_upper}
\nu(x) \le C\frac{\psi^*(|x|^{-1})}{|x|^d} , \quad x \in \bR^d\setminus\{0\}.
\end{equation}
\end{thm}

To obtain the estimates of density $p_t(x)$ and $\nu(x)$ we need assumptions on growth of $\psi$ near infinity
(see \cite[Section 3]{BGR}, \cite[Section 2]{KSV7} and  \cite[(2.7) and (2.20)]{Z}); 

\noindent
{\bf (H):}
there exist constants $a_1,a_2, r_0>0$ and  $\delta_1, \delta_2 \in (0,2)$ such that
\begin{equation}\label{e:S}
 a_1 s^{\delta_1} \le  \frac{\psi(st)}{\psi(t)} \le a_2 s^{\delta_2}  \qquad\textrm{ for all } s\ge 1 \textrm{ and } t \ge 1/r_0. 
\end{equation}

Recently, in \cite{BGR} Bogdan, Grzywny and Ryznar
 obtained an interesting equivalence on the upper and lower bounds of the densities. We state a partial result 
 relevant to our setting. 
\begin{thm}
[{\cite[Theorem 26]{BGR}}]\label{t:BGR}
 For an isotropic unimodal L\'evy process $X$ in $\bR^d$, the following are equivalent:
\begin{itemize}
\item[(i)]  {\bf (H)} hold for the characteristic function $\psi$ of $X$.
\item[(ii)] The transition density $p_t(x)$ of $X$ has following lower bound; for some $r_0\in(0,\infty)$ and a constant $c$,
$$
p_t(x) \ge c\frac{t\psi^*(|x|^{-1})}{|x|^d} , \quad 0< |x| < r_0, 0<t\psi^*(|x|^{-1})<1.
$$
\item[(iii)] The L\'evy density $\nu(x)$ of $X$ has following lower bound; for some $r_0\in(0,\infty)$ and a constant $c$,
\begin{equation}\label{nu_lower}
\nu(x) \ge c \frac{\psi^*(|x|^{-1})}{|x|^d}, \quad 0<|x| < r_0.
\end{equation}
\end{itemize}
\end{thm}

From
\eqref{e:almostinc}, \eqref{nu_upper} and Theorem \ref{t:BGR}, we conclude that our result cover isotropic unimodal L\'evy process
satisfying  {\bf (H)}. Recently the Harnack inequality and the H\"older estimate for harmonic functions with respect to isotropic unimodal L\'evy process was proved in \cite{G}.

A typical example of isotropic unimodal L\'evy process is a subordinate Brownian motion. If the characteristic exponent $\psi(r)$ is of the form $\phi(r^2)$ for some Bernstein function
\begin{equation*}
\phi(\lambda) = b\lambda + \int^{\infty}_0 (1-e^{-\lambda t}) \mu(dt),
\end{equation*}
where $b\ge 0 $ and $\mu$ is a measure on $(0,\infty)$ satisfying $\int_{(0,\infty)} (1 \wedge t ) \mu(dt) <\infty$, then the associated process $(X_t, t\ge 0)$ is a subordinate Brownian motion with L\'evy density
$$
\nu(x) = \int^{\infty}_0 (4\pi t)^{-d/2} e^{-|x|^2/(4t)} \mu(dt)
$$
and the infinitesimal generator of $X$ is $\phi(\Delta):=-\phi(-\Delta)$. 

One can find an extensive list of explicit Bernstein functions satisfying our assumption in \cite{SSV}.
Here are a few of them.

(1) $\phi(\lambda)=\lambda^{\alpha/2}$, $\alpha\in (0, 2)$
(symmetric $\alpha$-stable process);

(2) $\phi(\lambda)=(\lambda+m^{2/\alpha})^{\alpha/2}-m$, $\alpha\in (0, 2)$ and $m>0$
(relativistic $\alpha$-stable process);

(3) $\phi(\lambda)=
\lambda^{\alpha/2}+ \lambda^{\beta/2} $, $0<\beta<\alpha <2$
(mixed symmetric $\alpha$- and $\beta$-stable processes);

(4) $\phi(\lambda)=\lambda^{\alpha/2}(\log(1+\lambda))^{p}$,
 $\alpha\in (0, 2)$, $p\in [-\alpha/2, (2-\alpha)/2]$.

 \bigskip

{\bf Acknowledgement.} 
This paper is a part of the author's PhD thesis. He thanks Professor Panki Kim, his PhD thesis advisor, for his guidance and encouragement. 
The  author also wishes to thank Luis Silvestre who gave lectures on this regularity problem in the "Summer School 2012 on Nonlocal Operators" at the  Bielefeld University.

\end{doublespace}

\vskip 0.1truein

{\bf Jongchun Bae}

Department of Mathematical Sciences,

Seoul National University, Building 27, 1 Gwanak-ro, Gwanak-gu Seoul 151-747, Republic of Korea

E-mail: \texttt{bjc0204@snu.ac.kr}

\end{document}